\documentclass[12pt]{article}
\usepackage{graphics,latexsym,amssymb,amsmath,amscd}
\RequirePackage{graphics}
\usepackage[all,cmtip]{xy}
\def\[#1\]{\begin{eqnarray*}#1\end{eqnarray*}}

\def\tr{\hbox{\rm tr}\,}
\def\PU{\hbox{\bf PU}}
\def\SU{\hbox{\bf SU}}
\def\SL{\hbox{\bf SL}}
\def\phi{\varphi}

\input xy
\xyoption{all}
\newtheorem{thm}{Theorem}[section]
\newtheorem{dfn}[thm]{Definition}

\newtheorem{prop}[thm]{Proposition}

\newcommand{\Pf}{{\em Proof}. }
\newcommand{\EPf}{\hbox{}\hfill$\Box$\vspace{.5cm}}
\newcommand{\C}{{{\mathbb C}}}
\renewcommand{\P}{\mathbb {P}}
\newcommand{\R}{{{\mathbb R}}}
\newcommand{\Z}{{{\mathbb Z}}}

\def\tr{\hbox{\rm tr}\,}
\def\PU{\hbox{\bf PU}}
\def\SU{\hbox{\bf SU}}
\def\SL{\hbox{\bf SL}}
\def\PSL{\hbox{\bf PSL}}

\def\phi{\varphi}

\def\C{{\mathbb C}}
\def\R{{\mathbb R}}

\date{}
\title{Flag structures on real 3-manifolds}
\author{E. Falbel and
 J. M. Veloso}
\begin{document}
\maketitle
\newcommand{\D}{\mbox{$\cal D$}}
\newtheorem{df}{Definition}[section]
\newtheorem{te}{Theorem}[section]
\newtheorem{co}{Corollary}[section]
\newtheorem{po}{Proposition}[section]
\newtheorem{lem}{Lemma}[section]
\newcommand{\Ad}{\mbox{Ad}}
\newcommand{\ad}{\mbox{ad}}
\newcommand{\im}[1]{\mbox{\rm im\,$#1$}}
\newcommand{\bm}[1]{\mbox{\boldmath $#1$}}
\newcommand{\sime}{\mbox{sim}}
\begin{abstract}
We define flag structures on a real three manifold $M$ as the choice of two complex lines on the complexified tangent space at each point of $M$.   We suppose that the  plane field defined by the complex lines is a contact plane and construct an adapted connection on an appropriate principal bundle.  This includes path geometries and CR structures as special cases.  We prove that the null curvature models are given by totally real submanifolds in the flag space $\SL(3,\C)/B$, where $B$ is the subgroup of upper triangular matrices.   We also define a global invariant which is analogous to the Chern-Simons secondary class invariant for three manifolds with a Riemannian structure and to the Burns-Epstein invariant in the case of CR structures.  It turns out to be constant on homotopy classes of totally real immersions in flag space.

\end{abstract}

\section{Introduction}

Path geometries and CR structures on real three manifolds were studied by Elie Cartan in a long series of papers (see \cite{Ca,C} and \cite{B} for a beautiful account of this work).  Both geometries 
have models which are obtained through real forms of  a complex group.  More precisely, the 
group $\SL(3,\C)$ acts by projective transformations on both points in $ { \P({\C^3})}$ and its
lines viewed as ${ \P({\C^3}^*)}$.  The space of flags $F\subset { \P({\C^3})}\times \P({\C^3}^*)$ of lines containing points is described as the homegeneous space $\SL(3,\C)/B$ where $B$ is the subgroup of upper triangular matrices.  The path geometry of the flag space is defined by the
two projections onto points and lines in projective space.  Indeed, the kernels of the differentials of each projection define two complex lines in the tangent space of the flag space at each point.  It turns out that the planes generated in this way form a contact plane field.

The two real models appear as closed orbits of 
the two non-compact real forms $\SL(3,\R)$ and $\SU(2,1)$ in the space of flags.  
In this work we define a structure over a real manifold which interpolates between these two geometries.  Namely, the structure is a choice of two complex lines in the complexified tangent space at each point.  We call it a {\sl flag structure}.   Path geometry and CR  geometry correspond to a particular choice of 
lines adapted to the real structure of the real forms.   Flat structures modeled on  $\SU(2,1)$ are known as spherical CR structures and have been studied since Cartan.  But it is not known to what extent a 3-manifold may be equipped with such a structure (see \cite{BS1,S,DF} and their references).  A  flat path structure on a hyperbolic manifold was recently  constructed
in \cite{FS} where it is called a flag structure (as in \cite{Ba} where  flat path structures on Seifert manifolds arise  from  representations of surface groups into $\PSL(3,\R)$).  The use of configurations of flags associated to triangulations of 3-manifolds to obtain information on representation spaces is another theme related to this work (see \cite{BFG} and its references). 

Non-flat flag structures are abundant.  In fact, any 3-manifold has a real path structure and also a CR structure. 
We develop the equivalence problem for complex structures on a real three manifold in the first sections.   We obtain an adapted connection in an appropriate bundle over the manifold which has structure group $B$.   In  particular we characterize null curvature structures as those which can be locally embedded in the flag space $F$ and inherit the flag structure from the path structure of $F$. 
Real path structures and spherical CR structures of null curvature have a unique embedding up to translations by elements of $\SL(3,\C)$.

We define, using the adapted connection, a secondary invariant in the case the bundle has a global section.  The construction follows the same idea as  in \cite{BE} for the case of CR structures.  Computing the first variation of the invariant, we obtain, in particular, that it is invariant under deformations which are obtained from deformations of totally real immersions into flag space. 

Here is a more detailed account of each section.   
In section \ref{section:definition} we give the definition of flag structure and the classical examples of CR and 
real path geometry.  Those two geometries have been studied for a long time and very readable accounts are in \cite{J} for the CR case and \cite{IL} for real path geometry.  The basic examples of embedded totally real submanifolds in flag space are described in \ref{subsection:totallyreal}. Two families of homogeneous examples which are neither CR nor real path structures are described in the end of the section.

In section \ref{section:coframes} we define a coframe bundle $Y$ over a  $\C^*$-bundle $E$ over the real manifold adapted to the structure.  The $\C^*$-bundle is the set of contact forms at each point.  It is a complex line bundle with its zero section removed.  The coframe bundle $Y$ is then constructed over this bundle and the final descrition is given in Proposition \ref{proposition:coframe}.  As a principal bundle over the real manifold, $Y$ has structure group the group $B$ of upper triangular matrices (up to a finite cover).

Section \ref{section:connectionforms} is the technical core of the paper.  We construct the connection forms and curvature forms which will be intepreted as a Cartan connection in the following section.  The construction follows Cartan's technique
(see \cite{C,CM} and \cite{J}) but it is engineered to include  at the same time path geometries and CR structures.   The definition of Cartan connection we give (Definition \ref{definition:Cartan} in section \ref{section:Cartan}) is slightly more general than usual.  We don't impose that the tangent space of the fiber bundle be isomorphic to the Lie algebra.  This allows more flexibility and we are able to prove
Theorem \ref{theorem:Cartan} which puts together the construction in section \ref{section:connectionforms} into a Lie algebra valued connection form.    The characterization of null curvature structures is done in Theorem \ref{theorem:null} which identifies them locally as totally real embeddings into flag space. Finally we prove a rigidity theorem (Theorem \ref{theorem:rigidity}) which states that the only  CR or real path structures which admit local embeddings into the flag space are the flat ones.

In section \ref{section:global} we introduce our $\R$-valued global invariant for structures in the case the bundle $Y$ is trivial.  This follows Chern-Simons construction and coincides with Burns-Epstein invariant in the case of CR structures.  We do not attempt here to 
extend the construction to more refined versions as in \cite{CL} or \cite{BHR} for the CR case.  Also missing in this work is the application to second order differential 
equations which we plan to develop in a sequel of this paper.  The first variation formula (Proposition \ref{proposition:variation}) shows that
the critical points of the invariant occur at zero curvature structures.  In particular, using Theorem \ref{theorem:null}, we thus prove that totally real immersions are critical points with respect to the invariant. 
On the other hand, it follows from Gromov's h-principle techniques (Theorem 1.4 pg. 245 in \cite{Fo}) that every real 3-manifold admits a totally real immersion in $\C^3$ and therefore into the flag space (for the number of isotopy classes of totally real immersions see \cite{Bo}).  One can obtain in this way a set of numbers associated to a compact 3-manifold  corresponding to the values of the global invariant on  the space of totally real immersions with trivial bundle $Y$ up to homotopy. 

Section \ref{section:pseudo} concerns a natural reduction of the flag structure (which we call pseudo-flag structure).  Namely, when one choses a contact form, the structure group can be reduced and one can construct a bundle $X$ with an adapted connection.  In the case of CR structures this reduction (called pseudo-Hermitian structures) was throrougly analyzed in \cite{W} and we obtain, analogously, a particular embedding $X\rightarrow Y$ which allows 
one to describe our global invariant in terms of easier local data of the reduced bundle (see Proposition \ref{proposition:invariant}).   

In the last section we explicitly compute examples of homogeneous structures on $\SU(2)$.  The curvatures for a family of them
are non-zero and therefore, applying Theorem \ref{theorem:null}, they cannot be obtained as totally real manifolds in flag space. 
It would be interesting to understand if one can embed a flag structure in higher dimensional flag spaces.

We thank the R\'eseau Franco-Br\'esilien en Math\'ematiques for financial support during the preparation of this work.

\section{Flag structures in dimension 3}\label{section:definition}

Le $M$ be a real three dimensional manifold and $T^\C=TM\otimes \C$ be its complexified tangent bundle.

\begin{dfn}A flag structure on $M$ is a choice of two sub-bundles $T^1$ and $T^2$ in $T^\C$ such that 
$T^1 \cap T^2=\{ 0\}$ and such that $T^1 \oplus T^2$ is a contact distribution.
\end{dfn}

The condition that $T^1 \oplus T^2$ be a contact distribution means that, locally,  there exists a one form $\theta\in {T^*M}\otimes \C$ such that
$\ker \theta= T^1\oplus T^2$ and $d\theta \wedge \theta$ is never zero.

This definition contains two special cases, namely,
\begin{itemize}
 \item CR structures, which arise when $T^2=\overline {T^1}$.
 \item Path geometries, which are defined when $T^1$ and $T^2$ are complexifications of real one dimensional sub-bundles of $TM$.
\end{itemize}

Path geometries are treated in detail in section 8.6 of \cite{IL} and in \cite{BGH} where the relation to second order differential equations is also explained.   
CR structures in three dimensional manifolds were studied by Cartan (\cite{C}, see also \cite{J})  and his solution of the equivalence problem in dimension 3 was generalized to higher dimensions in \cite{CM}.  The goal in this paper is to use the same formalism 
for both geometries.  They appear as real forms of a complex path geometry.  Indeed, one can define a complex path geometry 
on complex manifolds of dimension three whose zero curvature model is the homogeneous flag space $\SL(3,\C)/B$ (where $B$ is the Borel subgroup of upper triangular matrices).
The group $\SL(3,\C)$ acts on this space transitively and the  real geometries associated to it correspond to the two closed orbits of the two non-compact  real forms $\SL(3,\R)$ and $\SU(2,1)$.

\subsection{The flag space }

The model cases of a flag structure arise when we consider certain  real three dimensional submanifolds of  the space $F$ of complete flags (that is, lines and planes containing them) in $\C^3$.  The group $\SL(3,\C)$ acts on the space of flags with isotropy $B$, the Borel group of upper triangular matrices.  We can describe therefore  the space of flags as the homogeneous space $F=\SL(3,\C)/B$.   The space of flags is equipped with two projections. One projects the line of a flag  into $\P(\C^3)$  on one hand and, on the other, the plane into $\P({{\C^3}^*})$ viewed as a kernel of a linear form.
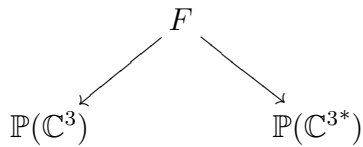
\begin{center}
$ 
\xymatrix{
      & F \ar@{>}[rd] \ar@{>}[ld]  & \\
       \P(\C^3)  &   & \P({\C^3}^*)\\
       }
$
\end{center}
The two projections correspond to the projections into $\P(\C^3)=\SL(3,\C)/P_1$ and $\P({\C^3}^*)=\SL(3,\C)/P_2$ where $P_1$ and $P_2$ are two different parabolic subgroups which fix, respectively, the line and the plane of the flag fixed by the Borel subgroup.. 

The projections define a complex contact distribution on the tangent space of $F$ generated by the tangent spaces to each of the fibers.  One can embed $F$ into $\P(\C^3)\times \P({\C^3}^*)$ as the set of pairs $(z,l)\in \P(\C^3)\times \P({\C^3}^*)$ satisfying the incidence relation $l(z)=0$.

Another description of the contact distribution can be given using  explicitly the Lie algebra structure of  $\SL(3,\C)$.
 Indeed, the Lie algebra of $\SL(3,\C)$ decomposes in the following direct sum of vector subspaces:
$$
 {\mathfrak {sl}}(3,\C)= {\mathfrak {g}}^{-2}\oplus {\mathfrak {g}}^{-1}\oplus {\mathfrak {g}}^{0}\oplus{\mathfrak {g}}^{1}\oplus {\mathfrak {g}}^{2}.
$$
That is  the graded decomposition of  ${\mathfrak {sl}}(3,\C)$ where ${\mathfrak b}= {\mathfrak {g}}^{0}\oplus{\mathfrak {g}}^{1}\oplus {\mathfrak {g}}^{2}$ corresponds to upper triangular matrices.  The tangent space of $ \SL(3,\C)/B$ at $[B]$ is identified to
$$
{\mathfrak {sl}}(3,\C)/{\mathfrak b}= {\mathfrak {g}}^{-2}\oplus {\mathfrak {g}}^{-1}.
$$ 
Moreover the choice of the reference flag whose isotropy is $B$ defines a decomposition of  ${\mathfrak {g}}^{-1}={\mathfrak {t_1}}^{}\oplus {\mathfrak {t_2}}^{}$, with $ {\mathfrak {g}}^{-2}=[{\mathfrak {t_1}},{\mathfrak {t_2}}]$, corresponding to the two parabolic subgroups with Lie algebras ${\mathfrak {p_1}}= {\mathfrak {t_1}}^{}\oplus {\mathfrak {g}}^{0}\oplus{\mathfrak {g}}^{1}\oplus {\mathfrak {g}}^{2}$ and ${\mathfrak {p_2}}= {\mathfrak {t_2}}^{}\oplus {\mathfrak {g}}^{0}\oplus{\mathfrak {g}}^{1}\oplus {\mathfrak {g}}^{2}$.

One propagates the decomposition of the tangent space at $[B]$ to the whole  flag space by the action of $\SL(3,\C)$ to obtain the two field of complex lines ${\mathfrak T}^1$ and ${\mathfrak T}^2$.


\subsection{Real submanifolds in flag space}\label{subsection:totallyreal}

In this section we assume $M$ is a three dimensional real manifold.  Let $\phi : M\rightarrow  \SL(3,\C)/B$ be an embedding and $\phi_* : TM\rightarrow T \SL(3,\C)/B$ be its differential.
One can extend this map to  $\phi_* : TM^\C\rightarrow T \SL(3,\C)/B$ by $\phi_*(iv)=J\phi_*(v)$ where $J$ is the complex structure on the tangent space of the complex manifold $\SL(3,\C)/B$ (which is just multiplication by $i$ in matrix coordinates). 

\begin{dfn}\label{definition:totallyreal}
An embedding $\phi : M\rightarrow  \SL(3,\C)/B$ is totally real if, for every $p\in M$, 
$$\phi_* : T_pM^\C\rightarrow  T_pSL(3,\C)/B$$ is an isomorphism.
\end{dfn}

One can define then the spaces $T^1, T^2\subset TM^\C$ which correspond to ${\mathfrak T}^1$ and ${\mathfrak T}^2$ on the flag space $\SL(3,\C)/B$.   This defines a flag structure on any real 3-manifold $M$ with a totally real embedding into $\SL(3,\C)/B$.

There are two fundamental examples.  They are both described as the unique closed orbit in the flag space by the action of a non-compact real form of $\SL(3,\C)$.

\subsubsection{Spherical CR geometry and $\SU(2,1)$}

Spherical CR geometry is modeled on the sphere ${\mathbb S}^3$ equipped with a natural ${\PU}(2,1)$ action.  Consider the group $\mathrm{U}(2,1)$ preserving the Hermitian form
$\langle z,w \rangle = w^*Jz$ defined on ${\C}^{3}$ by
the matrix
$$
J=\left ( \begin{array}{ccc}

                        0      &  0    &       1 \\

                        0       &  1    &       0\\

                        1       &  0    &       0

                \end{array} \right )
$$
and the following cones in ${\C}^{3}$:
$$
        V_0 = \left\{ z\in {\C}^{3}-\{0\}\ \ :\ \
 \langle z,z\rangle = 0 \ \right\},
$$
$$
        V_-   = \left\{ z\in {\C}^{3}\ \ :\ \ \langle z,z\rangle < 0
\ \right\}.
$$
Let $\pi :{\C}^{3}\setminus\{ 0\} \rightarrow  \P({\C^3})$ be the
canonical projection.  Then
${\mathbb H}_{\C}^{2} = \pi(V_-)$ is the complex hyperbolic space and its boundary is
$$
{\mathbb S}^3= \pi(V_0)=\{ [x,y,z]\in \mathbb{\C P}^{2}\ |\ x\bar z+ |y|^2+z\bar x=0\ \}.
$$
The group of biholomorphic transformations of ${\mathbb H}_{\C}^{2}$ is then
$\mathrm{\PU}(2,1)$, the projectivization of $\mathrm{U}(2,1)$.  
Observe that this group also acts on ${\mathbb S}^3$.

An element $x\in {\mathbb S}^3$ gives rise to a flag in $F$ where the line corresponds to the unique complex line tangent to ${\mathbb S}^3$ at $x$.
More explicitly, 
we let 
\begin{eqnarray*}
 \phi_{CR} \: : \: {\mathbb S}^3 & \to & F\\
 x & \mapsto &  \left(x, \langle\, . \,,\,  x \,\rangle \right)
\end{eqnarray*}

\subsubsection{Real path geometry and $\SL(3,\R)$}

Flat path geometry   is the geometry of real flags in $\R^3$.  That is the geometry of the space of all couples $[p,l]$ where $p\in \R P^2$ and $l$ is a real projective line
containing $p$.  The space of flags is identified to the quotient
$$
\SL(3,\R)/B_\R
$$
where $B_\R$ is the Borel group of all real upper triangular matrices.  The inclusion 
$$\phi_\R : \SL(3,\R)/B_\R\rightarrow \SL(3,\C)/B
$$
is clearly a totally real embedding. 

\subsubsection{Homogeneous immersions of $\SU(2)$ }

There exists two families:
\begin{enumerate}
\item
Let  $x_0=[1,0,0]\in  \P(\C^3)$ and define for a fixed $l\in \P({\C^3}^*)$ such that $l(x_0)=0$ the embedding.  Define the family of embeddings:

\begin{eqnarray*}
 \phi_{l} \: : \: { \SU(2)} & \to & F\\
 g & \mapsto &  \left(gx_0, g^*l \right)
\end{eqnarray*}
where $g\in \SU(2)\subset \SU(2,1)$.  Each choice of $l\in \P({\C^3}^*)$ defines a closed orbit in the space of flags by the action of $\SU(2)$ and therefore this family has $ \P({\C^2})$ as parameter space.  It can be seen as a deformation of the spherical CR structure on the sphere to a family of $\SU(2)$ invariant flag structures. The CR embedding is obtained when 
$l= \langle\, . \,,\,  x_0 \,\rangle$.  

\item  The other family arises from a deformation of the real flag space.
Let 
\begin{eqnarray*}
 \phi_{l} \: : \: {\mathbb{SO}}(3) & \to & F\\
 g & \mapsto &  \left(gx_0, g^*l \right)
\end{eqnarray*}
where $g\in {\mathbb{SO}}(3)\subset \SL(3,\R)$. 
\end{enumerate}
Again, each choice of $l\in \P({\C^3}^*)$ defines a closed orbit in the space of flags by the action of ${\mathbb{SO}}(3)$ and therefore this family has $\P(\C^2)$ as parameter space.  But there exists an isotropy. Namely, the orbit of $x_0$ is $\P(\R^3)$ and has isotropy ${\mathbb{O}}(2)$.  The deformation space in this case is the quotient ${\mathbb{O}}(2)\setminus  \P({\C^2})$ which is a segment.  The structures on $\SU(2)$ are obtained considering the twofold cover of ${\mathbb{SO}}(3)$.

\section{The $\C^*$-bundle of contact forms and an adapted coframe bundle}\label{section:coframes}

We start the construction of a canonical $\C^*$-bundle over a real three manifold equipped with a flag structure. 

We consider the forms $\theta$ on $T^\C$ such that $\ker \theta=T^1\oplus T^2$.  
Define ${E}$ to be the $\C^*$-bundle of all
such forms.   This bundle is trivial if and only if there exists a globally defined 
non-vanishing form $\theta$.

 On $E$ we define the tautological
form $\omega$.  That is $\omega_\theta=\pi^*(\theta)$ where
$\pi: {E}\rightarrow M$ is the natural projection.  

Fixing a form $\theta$ we next define forms $\theta^1$ and $\theta^2$ on $T^\C$ satisfying
$$
\theta^1(T^1)\neq 0 \ Ê \  
{\mbox{and}}\ 
\theta^2(T^2)\neq 0 .
$$
$$
\ker \theta^1=T^2 \ Ê \  {\mbox{and}}\ 
 \  \ker \theta^2=T^1.
$$
Fixing one choice, all others are given by $\theta'^{i} = a^i \theta^i + v^i\theta$.

We consider the tautological forms defined
by the forms above over the line bundle $E$.  That is, for each 
$\theta\in E$  we let $\omega^i_\theta= \pi^*(\theta^i)$.
At each point $\theta\in E$ we have the family of forms defined over $TE_{\theta}$
$$
        \omega' = \omega
$$
$$
        \omega'^{1} = a^1 \omega^1 + v^1\omega
$$
$$
        \omega'^{2} = a^2 \omega^2 + v^2\omega
$$

 We may, moreover, suppose that 
 $$
 d\theta= \theta^1\wedge \theta^2 \ \ \mbox{modulo} \ \theta
 $$ 
 and therefore
  $$
 d\omega= \omega^1\wedge \omega^2 \ \ \mbox{modulo} \ \omega .
 $$ 
  This imposes that $a^1a^2=1$.

  Those forms
vanish on vertical vectors, that is, vectors in the kernel
of the map $TE\rightarrow TM$. In order to define non-horizontal 1-forms 
we let $\theta$ be a section of $E$ over $M$ and introduce the coordinate
$\lambda\in \C^*$ in $E$.
By abuse of notation,
let $\theta$ denote the tautological form on the section $\theta$.
Therefore the tautological form $\omega$ over $E$ is 
$$
\omega_{\lambda}=\lambda \theta.
$$
Differentiating this formula we obtain
\begin{equation}  \label{domega}
                d\omega = \omega\wedge\phi + 
\omega^1\wedge\omega^{2}  
\end{equation}
where $\phi= -\frac{d\lambda}{\lambda}$ 
 modulo $\omega, \omega^1, 
\omega^{2}$.

 Observe that
$
\frac{d\lambda}{\lambda}
$ is a form intrinsically defined on $E$ up to horizontal
forms (the minus sign is just a matter of conventions).

For a different choice of forms satisfying the equation we write \ref{domega} as 
$$\begin{array}{rcl}
d\omega&=&\omega'\wedge\phi'+\omega'^{1}\wedge \omega'^{2}=\omega\wedge \phi'+(a^1 \omega^1 + v^1\omega)\wedge (a^2\omega^{2} + v^{2}\omega)\\
&=&
\omega\wedge (\phi'-a^1v^{2}\omega^1 +a^2v^1\omega^{2})+ \omega^1\wedge\omega^{2}
\end{array}
$$
it follows $\phi'=\phi +a^1v^{2}\omega^1 -a^2v^1\omega^{2} +s\omega$, with $s\in \C$.

We obtain in this way a coframe bundle  over $E$:
$$
        \omega' = \omega
$$
$$
        \omega'^{1} = a^1\omega^{1} + v^1 \omega
$$
$$
        \omega'^{2} = a^2\omega^{2} + v^2 \omega
$$
$$
        \phi'= \phi +a^1v^{2}\omega^1 -a^2v^1\omega^{2} +s\omega
$$
$v^1, v^2,s \in \C$  and $a^1,a^2\in \C^*$ such that $a^1a^2=1$.

\begin{dfn}
We denote by $Y$ the  coframe bundle
 $Y\rightarrow E$ given by the set of 1-forms  $\omega, \omega^{1},
  \omega^{2}, \phi$.
  Two coframes
are related by
$$
(\omega', \omega'^{1},
  \omega'^{2}, \phi')=
(\omega, \omega^{1},
  \omega^{2}, \phi)
 \left ( \begin{array}{cccc}

                        1       &    v^1                &       v^2 & s \\

                        0 &     a^1 & 0            &       a^1 v^2  \\
                        0   &  0 & a^2             &       -a^2v^{ 1} \\
                        0      &  0 & 0     &       1

                \end{array} \right )
$$
where  and $s, v^1, v^2 \in \C$ and $a^1,a^2\in \C^*$ satisfy $a^1a^2=1$.
\end{dfn}

The bundle $Y$ can also be fibered over the manifold $M$.  In order to describe the bundle $Y$ as a principal fiber bundle over
 $M$ observe that choosing a local section $\theta$ of $E$ and forms $\theta^1$ and $\theta^2$ on $M$ such that
 $d\theta=\theta^1\wedge \theta^2$ one can write a trivialization of the fiber 

$$
      \omega=  \lambda \theta
$$
$$
      \omega^1=   a^1\theta^{1} + v^1\lambda\theta
$$
$$
       \omega^2= a^2\theta^{2} + v^2\lambda\theta
$$
$$
       \phi= -\frac{d\lambda}{\lambda} +a^1v^{2}\theta^1 -a^2v^1\theta^{2} +s\theta,
         $$
where $v^1, v^2,s \in \C$  and $a^1,a^2\in \C^*$ such that $a^1a^2=\lambda$.  Here the coframe $\omega, \omega^1,\omega^2,\phi$ is seen as 
composed of tautological forms.

The group $H$ acting on the right of this bundle 
 is 
$$
H=\left\{ \left ( \begin{array}{cccc}

                        \lambda   	&    	v^1 \lambda      	&       v^2\lambda 	& 		s \\

                       0				 &     a^1 			&		 0            &       a^1 v^2\\
                       0   				&  	0 				&		 a^2        &        -a^2v^{ 1}\\
                        0		       &  		0              &     0   			&       1

                \end{array} \right )
\mbox{
where   $s, v^1, v^2 \in \C$ and $a^1,a^2\in \C^*$ satisfy $a^1a^2=\lambda$
}
\right\}.
$$

Consider the homomorphism from the Borel group $B\subset \SL(3,\C)$ of upper triangular matrices into $H$
$$
j: B\rightarrow H
$$
given by

$$
\left( \begin{array}{ccc}

                        a  	&    	c      	&       e \\

                       0				 &     \frac{1}{ab}		&		d\\
                       0   				&  	0 				&		b
                      
       \end{array} \right ) 
       \longrightarrow  
       \left ( \begin{array}{cccc}

                        \frac{a}{b}  		&    	-2a^2d      	&       	\frac{c}{b}		& 		\frac{4e}{b}-2acd \\

                       0				&     a^2b 		&		0            		&      		abc\\
                       0   			&  	0 		&		\frac{1}{ab^2}    &        	\frac{2d}{b}\\
                       0		         	&  	0              	&       	0   			&       	1

                \end{array} \right )
$$
One verifies that the homomorphism is surjective and its kernel is isomorphic to $\Z/3\Z$ formed by diagonal matrices, 
so  that $H$ is isomorphic to the Borel group of projected upper triangular matrices in $\PSL(3,\C)$.

\begin{prop}  The bundle $Y\rightarrow M$ is a principal bundle with structure group $B/\C^*$ where $B$ is the Borel group of upper triangular matrices.
\end{prop}\label{proposition:coframe}

\section{Construction of connection forms on the bundle $Y$}\label{section:connectionforms}

 The goal of this section is to  obtain canonical forms defined on the coframe bundle $Y\rightarrow E$.  They will
 give rise to a connection on $Y$ as explained in the next section.  The connection will be a form on $Y$ with values in $   {\mathfrak {sl}}(3,\C)$ and will not be a Cartan connection as it is usually defined but a slight generalization of it.

 A local section of the complexified coframe bundle over $M$ may be given by three forms 
 $$
 \theta, \ \ \theta^1, \ \ \theta^2
 $$
 satisfying 
 $d\theta=\theta^1\wedge \theta^2$, with $\ker \theta^1=T^1$ and  $\ker \theta^2=T^2$.
 They give coordinates on the complexified cotangent bundle over $E$ and, furthermore,  we may describe next the tautological forms
 defined over that bundle in the previous section.
 
  
  At the point $\lambda\theta\in E$, the coframes of $Y$ are parametrized as follows:
 $$
 \omega= \lambda \theta
 $$
 $$
 \omega^i=a^i\theta^i + v^i\lambda\theta
 $$
 with $a^1a^2=\lambda$ and
 $$
 d\omega =\omega^1\wedge \omega^2+ \omega\wedge \phi
 $$
 where $\phi = -\frac{d\lambda}{\lambda}  \  \mod{\omega^1,\omega^2,\omega}$.

 We consider now these forms as tautological forms on the bundle $T^*Y$.
Differentiating $\omega^1,\omega^2$ we obtain
$$
d\omega^1= da^1\wedge \theta^1 + dv^1\wedge \lambda\theta +a^1 d\theta^1+v^1d\,\lambda\theta
$$
$$
d\omega^2= da^2\wedge \theta^2+ dv^2\wedge \lambda\theta +a^1 d\theta^2+v^2d\,\lambda\theta.
$$
Observing now that 
$$
\frac{d\lambda}{\lambda}= \frac{da^1}{a^1}+\frac{da^2}{a^2}
$$
we can write (modulo $\omega^1,\omega^2,\omega$)
$$
d\omega^1\equiv \frac{d\lambda}{2\lambda} \wedge \omega^1+ \frac{1}{2}\left(\frac{da^1}{a^1}-\frac{da^2}{a^2}\right)\wedge \omega^1 
+{dv^1}\wedge\omega
$$
$$
d\omega^2 \equiv \frac{d\lambda}{2\lambda}\wedge  \omega^2- \frac{1}{2}\left(\frac{da^1}{a^1}-\frac{da^2}{a^2}\right)\wedge \omega^2
+{dv^2}\wedge\omega
$$
Now, distributing the missing terms in $\omega^1,\omega^2,\omega$ in the last two terms of $d\omega^1$ and $d\omega^2$ and anti-symmetrizing
we obtain the following:

\begin{lem}
There exists linearly independent forms $\omega^1_1, \phi^1,\phi^2$ defined on $T^*Y$
such that
\begin{equation}\label{domegai}
d\omega^1=\frac{1}{2} \omega^1\wedge \phi +
\omega^1\wedge\omega^1_1 +\omega\wedge \phi^1\  \mbox{and}\ \ 
d\omega^2=\frac{1}{2} \omega^2\wedge \phi -
\omega^2\wedge\omega^1_1 +\omega\wedge \phi^2
\end{equation}
\end{lem}
Let  ${{\omega'}^1_1}$, ${{\phi'}^1}$ and ${\phi'^2}$ be other forms satisfying
equations \ref{domegai}. Taking the difference
we obtain for $i=1, 2$
$$
0=\omega^1\wedge(\omega^1_1-\omega'^1_1) +\omega\wedge (\phi^1-\phi'^1)
$$
and
$$
0=-\omega^2\wedge(\omega^1_1-\omega'^1_1) +\omega\wedge (\phi^2-\phi'^2)
$$
Therefore
$$
\omega^1_1-\omega'^1_1 = a\omega
$$
$$
\phi^1-\phi'^1 = a\omega^1 +b^1\omega
$$
$$
\phi^2-\phi'^2 = -a\omega^2 +b^2\omega
$$.  
\begin{lem}
There exists a 1-form $\psi$ such that
\begin{equation}\label{dphi}
d\phi = 
\omega ^1\wedge\phi^{2} -
\omega^{2}\wedge \phi^1 +\omega \wedge \psi
\end{equation}
\end{lem}
\Pf
Differentiating equation 
$$
d\omega =\omega^1\wedge \omega^2+ \omega\wedge \phi
$$
and using equations \ref{domegai}
we obtain
$$
\omega\wedge \left ( d\phi +\omega^{2}\wedge \phi^1 -
\omega ^1 \wedge \phi^2\right ) =0
$$which implies that there exists a form $\psi$ as claimed.
\EPf

If other forms  $\psi'$, ${{\phi'}^1}$ and ${\phi'^2}$ satisfy \ref{dphi}  then
$$
0=\omega ^1\wedge(\phi^{2}-\phi'^{2}) -
\omega^{2}\wedge( \phi^1-\phi'^1) +\omega \wedge( \psi-\psi')$$ 
and, therefore, using  $\phi^1-\phi'^1 = a\omega^1 +b^1\omega$ and
$\phi^2-\phi'^2 = -a\omega^2 +b^2\omega$ we obtain

$$
\psi-\psi'=b^2 \omega^1-b^1 \omega^{2}+c\omega .
$$
Our next goal is  to fix $\omega^1_1$.  For that sake we differentiate equations \ref{domegai}.
Equation $dd\omega^1=0$ gives
$$
\omega^1\wedge \left (-d\omega^1_1+\frac{3}{2} \omega^2\wedge \phi^1\right ) +
 \omega \wedge \left ( -d\phi^1 +\phi^1\wedge \omega^1_1+\frac{1}{2}\omega^1\wedge \psi-\frac{1}{2}\phi^1\wedge \phi\right )=0.
$$
Analogously, equation $dd\omega^2=0$ gives
$$
\omega^2\wedge \left (d\omega^1_1-\frac{3}{2} \omega^1\wedge \phi^2\right ) +
 \omega \wedge \left ( -d\phi^2 -\phi^2\wedge \omega^1_1+\frac{1}{2}\omega^2\wedge \psi-\frac{1}{2}\phi^2\wedge \phi\right )=0.
$$
Defining
$$
\Omega^1_1=d\omega^1_1-\frac{3}{2} \omega^2\wedge \phi^1 -\frac{3}{2} \omega^1\wedge \phi^2
$$
$$
\Phi^1= d\phi^1 -\phi^1\wedge \omega^1_1-\frac{1}{2}\omega^1\wedge \psi-\frac{1}{2}\phi\wedge \phi^1 
$$
$$
\Phi^2=d\phi^2 +\phi^2\wedge \omega^1_1-\frac{1}{2}\omega^2\wedge \psi-\frac{1}{2}\phi\wedge \phi^2
$$ 
So the equations can be written
$$
\omega^1\wedge \Omega^1_1+\omega\wedge \Phi^1=0
$$
$$
\omega^2\wedge \Omega^1_1-\omega\wedge \Phi^2=0
$$
The first equation implies that $\Omega^1_1=\omega\wedge \lambda^1_1+\omega^1\wedge \mu$, where $\lambda^1_1$ and $\mu$ are 1-forms
so that $\mu$ has no terms in $\omega$.  From the second equation we obtain that $\mu=S_1\omega^1+S_2\omega^2$. Substituting in the second we see that $S_1=0$.  We write therefore
$$
\Omega^1_1=\omega\wedge \lambda^1_1+S_{2}\omega^1\wedge \omega^2.
$$

\begin{lem}
There exists a unique form $\omega^1_1$ such that $\Omega^1_1=\omega\wedge \lambda^1_1$ (that is $S_{2}=0$).
\end{lem}
\Pf
Computing the difference 
$$
 \Omega^1_1- \Omega'^1_1=d(\omega^1_1-\omega'^1_1)-\frac{3}{2} \omega^2\wedge (\phi^1-\phi'^1) -\frac{3}{2} \omega^1\wedge (\phi^2-\phi'^2)
 $$
 and obtain
 $$
  \Omega^1_1- \Omega'^1_1=4a\,\omega^1\wedge \omega^2\ \ \mod \omega.
  $$
One can fix therefore $a$ so that $\Omega^1_1$ satisfies the assertion of the lemma.
\EPf
\subsubsection{} As $\omega^1_1$ is fixed we still have the following ambiguities:
$$
\phi^1-\phi'^1 = b^1\omega
$$ 
$$
\phi^2-\phi'^2 = b^2\omega
$$
$$
\psi-\psi'=b^2 \omega^1-b^1 \omega^{2}+c\omega .
$$
\begin{lem}
There exists unique forms $\phi^1$ and $\phi^2$ such that $\Omega^1_1$ does not contain terms  $\omega^i\wedge \omega$, $i=1,2$.
\end{lem}
\Pf
From the definition of $\Omega^1_1=d\omega^1_1-\frac{3}{2} \omega^2\wedge \phi^1 -\frac{3}{2} \omega^1\wedge \phi^2$ 
we obtain that
$$
\Omega^1_1-\Omega'^1_1=-\frac{3}{2}b^1 \omega^2\wedge \omega -\frac{3}{2}b^2 \omega^1\wedge \omega.
$$
So we can choose $b^1$ and $b^2$ as claimed.
\EPf

Observe that we fixed the 1-forms $\omega^1_1, \phi^1$ and $\phi^2$ so that
\begin{equation}  \label{domega11}
d\omega^1_1-\frac{3}{2} \omega^2\wedge \phi^1 -\frac{3}{2} \omega^1\wedge \phi^2=\omega\wedge \lambda^1_1
\end{equation}
with $\lambda^1_1\equiv 0 \mod \phi, \phi^1,\phi^2,\psi$.

\subsubsection{}
It remains to fix $\psi$.

\begin{lem}
There exists a unique 1-form $\psi$ such that $\Phi^1$ does not contain a term $\omega^1\wedge \omega$.
\end{lem}
\Pf
We compute using the definition
$\Phi^1-\Phi'^1= -\frac{1}{2}\omega^1\wedge (\psi-\psi')=-\frac{1}{2}c\omega^1\wedge \omega$
and we can choose a unique $c\in \C$ which proves the lemma.
\EPf 
 
 
 \subsection{Curvature forms}
 
 Curvature forms appear as differentials of connection forms and were used implicitly in the previous paragraphs to 
 fix the connection forms.
 
 In this section we obtain properties of the curvature forms which will be used in the following sections.
 Substituting $\Omega^1_1=\omega\wedge \lambda^1_1$ in equation $\omega^1\wedge \Omega^1_1+\omega \wedge \Phi^1=0$
 we obtain  that $\omega\wedge ( -\omega^1\wedge \lambda^1_1+ \Phi^1)=0$ and therefore 
 $\Phi^1=\omega^1\wedge \lambda^1_1+\omega\wedge \nu^1$ for a 1-form $\nu^1$.  Observe that 
 $\nu^1\equiv 0 \mod \omega^2, \phi, \phi^1,\phi^2, \psi$ in view of the last lemma.
 Analogously, one may write $\Phi^2=-\omega^2\wedge \lambda^1_1+\omega\wedge \nu^2$.
 
 \subsubsection{}
 Equation $d(d\phi)=0$ obtained differentiating \ref{dphi} can be simplified to 
 $$
0= \omega\wedge \left( -d\psi +\phi\wedge \psi + 2\phi^1\wedge \phi^2 -\omega^1\wedge \nu^2+\omega^2\wedge \nu^1\right).
 $$
 It follows that there exists a 1-form $\nu$ such that
 \begin{equation}\label{dpsi}
 d\psi-2\phi^1\wedge \phi^2-\phi\wedge \psi +\omega^1\wedge \nu^2-\omega^2\wedge \nu^1=\nu\wedge \omega.
 \end{equation}
 
 \subsubsection{}
 Equation $d(d\omega^1_1)=0$ obtained differentiating \ref{domega11} can be simplified to
 $$
 0=\omega^1\wedge \omega^2 \wedge 4\lambda^1_1+\omega\wedge 
 \left(-d\lambda^1_1+\frac{3}{2}\omega^1\wedge \nu^2+\frac{3}{2}\omega^2\wedge \nu^1+\phi\wedge \lambda^1_1\right)
 $$
 which implies, as $\lambda^1_1$ does not contain terms in $\omega$ that $\lambda^1_1=0$ and
  $\omega\wedge (\omega^1\wedge \nu^2+\omega^2\wedge \nu^1)=0$.  As $\nu^1$ does not have a term in $\omega^1$
  it follows from the last equation that $\nu^2$ does not have a term in $\omega^2$ and we conclude that
  $\nu^1=Q^1\omega^2$ and $\nu^2=Q^2\omega^1$ where we introduce functions $Q^1$ and $Q^2$.
  
  We have obtained the following equations:
  \begin{equation}\label{omega11}
\Omega^1_1:=d\omega^1_1-\frac{3}{2} \omega^2\wedge \phi^1 -\frac{3}{2} \omega^1\wedge \phi^2=0,
\end{equation}
 \begin{equation}\label{Phi1}
\Phi^1:= d\phi^1 -\phi^1\wedge \omega^1_1-\frac{1}{2}\omega^1\wedge \psi-\frac{1}{2}\phi\wedge \phi^1 =Q^1\omega\wedge \omega^2,
\end{equation}
 \begin{equation}\label{Phi2}
\Phi^2:=d\phi^2 +\phi^2\wedge \omega^1_1-\frac{1}{2}\omega^2\wedge \psi-\frac{1}{2}\phi\wedge \phi^2=Q^2\omega\wedge \omega^1.
\end{equation}

  \subsubsection{}
  Equation $d(d\phi^1)=0$ obtained differentiating $\Phi^1$ above is simplified to
  $$
  0=\frac{1}{2}\omega\wedge\omega^1\wedge \nu 
 - \omega\wedge \omega^2\wedge \left( dQ^1 +2Q^1\omega^1_1 -2Q^1\phi\right).
 $$
 It implies that 
 $$
 \nu= U_1\omega^1+ U_2\omega^2
 $$
and
$$
   dQ^1 +2Q^1\omega^1_1 -2Q^1\phi=S^1\omega -\frac{1}{2}U_2\omega^1 + T^1\omega^2,
   $$
   where we introduced functions $U_1, U_2, S^1$ and $T^1$.
   
   We obtain, substituting the expression for $\nu$ in equation \ref{dpsi}, the following expression
   \begin{equation}\label{dPsi}
\Psi := d\psi-2\phi^1\wedge \phi^2-\phi\wedge \psi +\omega^1\wedge \nu^2-\omega^2\wedge \nu^1=(U_1\omega^1+ U_2\omega^2)\wedge \omega.
 \end{equation}
 
 Therefore, as $\nu^1=Q^1\omega^1$ and $\nu^2=Q^2\omega^2$ we have
  \begin{equation}\label{dPsi'}
\Psi := d\psi-2\phi^1\wedge \phi^2-\phi\wedge \psi =(U_1\omega^1+ U_2\omega^2)\wedge \omega.
 \end{equation}

   \subsubsection{}
   Anagously, equation $d(d\phi^2)=0$ obtained differentiating $\Phi^2$ above is simplified to
 $$
  0= 
  \omega\wedge \omega^1\wedge \left( dQ^2 -2Q^2\omega^1_1 -2Q^2\phi +\frac{1}{2}U_1\omega^2\right).
 $$
 It implies that 
 $$
   dQ^2 -2Q^2\omega^1_1 -2Q^2\phi=S^2\omega -\frac{1}{2}U_1\omega^2 + T^2\omega^1,
   $$
   where we introduced new functions $S^2$ and $T^2$.

\subsubsection{} Finally, equation $d(d\psi)$ obtained from \ref{dpsi} simplifies to
$$
0=\omega\wedge \omega^1\left( dU_1-\frac{5}{2} U_1\phi -U_1\omega^1_1 +2Q^2\phi^1\right)
+\omega\wedge \omega^2\left( dU_2-\frac{5}{2} U_2\phi +U_2\omega^1_1 -2Q^1\phi^2\right)
 $$
 which implies that
 $$
 dU_1-\frac{5}{2} U_1\phi -U_1\omega^1_1 +2Q^2\phi^1=A\omega+B\omega^1+C\omega^2
 $$
 and
 $$
 dU_2-\frac{5}{2} U_2\phi +U_2\omega^1_1 -2Q^1\phi^2=D\omega+C\omega^1+E\omega^2.
 $$
\subsection{The CR case}

In order to make the construction of the bundle $Y$  compatible with the  CR bundle constructed 
in \cite{C,CM} one verifies first that $\omega^2=i\bar\omega^1$.   We have indeed
$$
d\omega=\omega\wedge \phi+i\omega^1\wedge\bar\omega^1.
$$
The form $\omega$ can be taken to be real so $\phi$ is also real.

From equations \ref{domegai} one has
$$
\omega^1_1+\bar \omega^1_1 =0 \ \ \ \phi^2=i\bar \phi^1.
$$
 From equation \ref{omega11}
we obtain
$$
d\omega^1_1-\frac{3}{2} i\bar \omega^1\wedge \phi^1 -\frac{3}{2} \omega^1\wedge\bar \phi^1=0,
$$
From equations \ref{Phi1},\ref{Phi2}
$$
d\phi^1 -\phi^1\wedge \omega^1_1-\frac{1}{2}\omega^1\wedge \psi-\frac{1}{2}\phi\wedge \phi^1 =Q^1\omega\wedge \omega^2,
$$
$$
d\phi^2 +\phi^2\wedge \omega^1_1-\frac{1}{2}\omega^2\wedge \psi-\frac{1}{2}\phi\wedge \phi^2=Q^2\omega\wedge \omega^1,
$$
we have $\psi=\bar \psi$ and $Q^1=\bar Q^2$.
From equation \ref{dPsi'} we obtain
$$
 d\psi-2\phi^1\wedge \bar \phi^1-\phi\wedge \psi -(U_1\omega^1+ U_2\bar \omega^1)\wedge \omega=0.
 $$
 From that equation we observe that $U_1=-i\bar U_2$ (in the CR literature $U_1$ is denoted $R_1$).

\section{The Cartan connection}\label{section:Cartan}

We consider the bundle $Y\rightarrow M$ as a principal bundle with structure group $H\subset \SL(3,\C)$, the Borel group of triangular matrices.
Observe that, although $M$ is a real manifold, each fiber is a complex space of dimension five which can be identified to the Borel 
subgroup $H$.    The real dimension of $Y$ is 13 which is $\dim_\R \SL(3,\C)-3$.  This dimension difference does not allow us to obtain a genuine Cartan connection but a slight generalization of the definition will be sufficient for our purposes. 

Recall that $X^*(y)=\frac{d }{dt}_{_{t=0}}{y e^{tX}}$ where $e^{tX}$ is the one parameter group generated
by  $X$.

\begin{dfn}\label{definition:Cartan}

A Cartan connection on $Y$ is a 1-form $\pi: TY \rightarrow {\mathfrak {sl}}(3,\C)$ satisfying:\newline
0. $\pi_p:T_pY\rightarrow  {\mathfrak {sl}}(3,\C)$ is injective.\newline
	1. If $X\in {\mathfrak h}$ and $X^*\in TY$ is the vertical vector field  canonically associated to $X$ then
	$\pi(X^*)=X.$\newline
	2. If $h\in H$ then  $(R_h)^* \pi=Ad_{h^{-1}}\pi$
\end{dfn}

Note that contrary to the usual definition of Cartan connection we don't impose that $\pi_p:T_pY\rightarrow  {\mathfrak {sl}}(3,\C)$ be an isomorphism as the dimensions are different.
 
We can represent the structure equations \ref{domega}, \ref{domegai}, \ref{dphi}, \ref{domega11} and \ref{dpsi}
  as a matrix equation whose entries are differential forms. The forms are disposed in the Lie algebra ${\mathfrak {sl}}(3,\C)$ as

$$\pi=\left ( \begin{array}{ccc}
                        -\frac{1}{2}\phi-\frac{1}{3}\omega^1_1     &   \phi^{ 2}   &   -\frac{1}{4}\psi  \\

                     \omega^1          &   \frac{2}{3}\omega^{1}_{1}    &    \frac{1}{2}\phi^1  \\

                      2\omega         &  2\omega^{2}   &      \frac{1}{2}\phi-\frac{1}{3}\omega^1_1
                \end{array} \right )
$$

It is a simple verification to show that 
\begin{equation}\label{dpi}
d\pi+\pi\wedge\pi=\Pi
\end{equation}
where

$$\Pi=\left ( \begin{array}{ccc}
                        0 &   -\Phi^{2}      &  -\frac{1}{4}\Psi   \\

                         0  &  0    &  \frac{1}{2}\Phi^1  \\

                           0  &   0  &    0
               \end{array} \right )
$$

with $
\Phi^1=Q^1\omega\wedge \omega^2,\ \ 
\Phi^2=Q^2\omega\wedge \omega^1 \ \mbox{and}\ \ \Psi =\left( U_1\omega^1+ U_2\omega^2\right) \wedge \omega.
$

\begin{thm}\label{theorem:Cartan}
The form $\pi$ is a Cartan connection on $Y\rightarrow M$.
\end{thm}

\Pf 
The action by $H$ can be replaced by the action of the Borel group of upper triangular matrices as described in a previous section.
The action on the right by an element $h\in H$ on a coframe $y\in Y$ is denoted by $R_h(y)$. As $y\in Y$ is a coframe of $E$, one may 
consider tautological lifts $y^{taut}$  on $Y$ defined by the elements of the coframe $y$.  The action on $Y$ lifts to an action on its tautological lifts as follows:
$$
{R_h}^*(y^{taut})={R_h (y)}^{taut}.
$$
We  compute  $Ad_{h^{-1}}\pi$ for an element 
$$
h=
\left( \begin{array}{ccc}

                        a  	&    	c      	&       e \\

                       0				 &     \frac{1}{ab}		&		d\\
                       0   				&  	0 				&		b
                      
       \end{array} \right )
$$
and verify that the  tautological forms $\omega,\omega^1,\omega^2,\phi$ on $Y$ (which appear as certain components of the connection) change according to 
the right action above.

It remains to show that the other components change similarly.
Now given $\omega,\omega^1,\omega^2,\phi$ tautological forms on $Y$ we defined unique forms $\omega^1_1,\phi^1,\phi^2,\psi$ such that the curvatures
$\Phi^1,\Phi^2$ and $\Psi$ had special properties.

We have 
$$
Ad_{h^{-1}}(d\pi-\pi\wedge\pi)=dAd_{h^{-1}}\pi-Ad_{h^{-1}}\pi\wedge Ad_{h^{-1}}\pi.
$$
Writing 
$$
\tilde{\pi}=Ad_{h^{-1}}\pi,
$$
we obtain $Ad_{h^{-1}}\Pi=\tilde \Pi$, where $\tilde \Pi=d\tilde\pi-\tilde\pi\wedge \tilde\pi$.  A computation shows that

\begin{align}
\label{adjunta}
\tilde\omega  &= \frac{a}{b}\, \omega\notag\\
\tilde{\omega}^1 &=a^2b\, \omega^1-2da^2\,\omega\notag\\
\tilde\omega^2 &=\frac{1}{ab^2}\,\omega^2+\frac{c}{b}\,\omega\notag\\
\tilde\phi &=\phi+abc\,\omega^1+2\frac{d}{b}\,\omega^2+(\frac{4e}{b}-2dac)\,\omega\\
\tilde\omega^1_1 &=\omega^1_1+\frac{3}{2}abc\,\omega^1-3\frac{d}{b}\,\omega^2-3dac\,\omega\notag\\
\tilde\phi^1 &=b^2a\, \phi^1+2dab\, \omega^1_1-bad\,\phi+2bae\, \omega^1-4d^2a\,\omega^2-4dae\, \omega\notag\\
\tilde\phi^2 &=\frac{1}{ba^2}\, \phi^2+\frac{c}{a}\, \omega^1_1+\frac{1}{2}ca\,\phi+bc^2\, 
\omega^1+(\frac{2e}{a^2b^2}-\frac{2cd}{ab})\,\omega^2+(\frac{2ce}{ab}-{2dc^2})\, \omega\notag\\
\tilde\psi &=\frac{b}{a}\,\psi+(\frac{4e}{a}-2bc)\,\phi +4bce\,\omega^1+(\frac{8de}{ab}-8{cd^2})\,
\omega^2+2cb^2\,\phi^1+4\frac{d}{a}\, \phi^2+4dbc\, \omega^1_1+(\frac{8e^2}{ab}-{8dce})\, \omega\notag
\end{align}

therefore it suffices 
to  verify that the new curvature forms  $\Phi'^1,\Phi'^2$ and $\Psi'$ obtained from $Ad_{h^{-1}}\Pi$ verify the same properties.
Indeed 
$$
Ad_{h^{-1}}\Pi=
\tilde \Pi=\left ( \begin{array}{ccc}
                        0 &  \frac{1}{a^2b}\Pi_{12}     & \frac{b}{a}\Pi_{13}+ \frac{d}{a}\Pi_{12}-cb^2\Pi_{23}\\

                         0  &  0    &  ab^2\Pi_{23}  \\

                           0  &   0  &    0
               \end{array} \right )
$$
%
%
and we see  that the new curvatures satisfy the same properties.
\EPf

In particular, one obtains that
$$
\tilde \Phi^1=2\tilde \Pi_{23}=2ab^2\Pi_{23}\Pi_{23}=ab^2\Phi^1
$$
and therefore
$$
\tilde Q^1\tilde \omega\wedge \tilde\omega^2=ab^2\, Q^1\omega\wedge\omega^2.
$$
But $\tilde Q^1\tilde \omega\wedge \tilde\omega^2=\tilde Q^1\frac{a}{b}\, \omega\wedge \frac{1}{ab^2}\omega^2$
and then
$$
\tilde Q^1=ab^5\, Q^1.
$$
Analogously,
from
 $$
\tilde \Phi^2=-\tilde \Pi_{12}=-\frac{1}{a^2b}\Pi_{12}\Pi_{23}=\frac{1}{a^2b}\Phi^2
$$
we obtain that
$$
\tilde Q^2=\frac{1}{a^5b}\, Q^2.
$$
These transformation properties imply that we can define two tensors on $Y$ which are invariant under $H$ and will give rise to two functions  on $M$.
Indeed
$$
 Q^1\,\omega^2\wedge  \omega\otimes { \omega}\otimes   e_1
$$
and
$$
 Q^2\,\omega^1\wedge  \omega\otimes { \omega} \otimes  e_2,
$$
where $e_1$ and $e_2$ are duals to $\omega^1$ and $\omega^2$ in the dual frame of the coframe bundle of $Y$ are easily seen to be $H$-invariant.

\subsection{Null curvature models}

 A local characterization of null curvature 
path geometries is given in the following theorem.
Recall definition \ref{definition:totallyreal} of 
a totally real embedding $\phi : M\rightarrow  \SL(3,\C)/B$ and 
its associated  flag structure  which correspond to $\phi_*(T^1)={\mathfrak {T_1}}$ and $\phi_*(T^2)={\mathfrak {T_2}}$ on the flag space $\SL(3,\C)/B$.

\begin{thm}\label{theorem:null}
A totally real embedding $\phi : N\rightarrow  \SL(3,\C)/B$ with induced  flag structure on $ TM^\C$ given by $T^1$ and $T^2$ as above
 is a contact path structure with adapted connection having null curvature.
Conversely a contact path structure whose adapted connection has zero curvature is locally equivalent to a totally  real embedding with induced path structure defined by $T^1$ and $T^2$ as above.
\end{thm}

Observe that null curvature does not define a unique flag structure on a real manifold but instead decides whether is can be embedded as a totally real submanifold in flag space. 
\vspace{.5cm}

\Pf  
The fact that $N$, equipped with the  two sub-bundles $T^1$ and $T^2$ is a contact path structure with adapted connection having null curvature
follows from the fact that the adapted principal bundle $Y$ associated to $N$ is identified to the restriction to $N$ of the bundle  $\SL(3,\C)\rightarrow \SL(3,\C)/H$.  The adapted connection on $Y$ is then the Maurer-Cartan form of $\SL(3,\C)$ restricted to this bundle and therefore has zero curvature.

Suppose now that $M$ has a contact path structure defined by sub-bundles $T^1$ and $T^2$ in $TN\otimes \C$. 
Let $\pi: TY \rightarrow {\mathfrak {sl}}(3,\C)$ be an adapted connection.  Suppose that $d\pi+\pi\wedge\pi=\Pi=0$ and let $\tilde \omega$ be the 
Maurer-Cartan form of the group $\SL(3,\C)$.  By Cartan's theorem (see theorem 1.6.10 in \cite{IL})  every $y\in Y$ is contained in an open neighborhood 
$U\subset Y$ where an immersion $f: U\rightarrow \SL(3,\C)$ is defined satisfying $\pi=f^*(\tilde \omega)$.  Moreover, any other immersion $\tilde f$ satisfying the same equation 
is related by a translation by an element $a\in \SL(3,\C)$ in the group, that is, $\tilde f= a f$.  

If $X\in {\mathfrak h}$ then 
$X=\pi(X^*)=f^*(\tilde \omega)(X^*)=\tilde\omega(f_*X^*)$.  Therefore $f_*X^*$ is tangent to the fibers of $\SL(3,\C)\rightarrow \SL(3,\C)/H$. We conclude that
$f: U\rightarrow \SL(3,\C)$ projects to an immersion $\bar f: V\rightarrow \SL(3,\C)/H$ where $V\subset M$.
The subspaces $\bar f_* (T^1)$ and $\bar f_* (T^2)$ are precisely the subspaces ${\mathfrak {T_1}}$ and ${\mathfrak {T_2}}$ restricted to $\bar f(V)$.

\EPf

The following theorem shows the rigidity of the real models in flag space.  It remains the possibility that
general CR structures or path structures might be deformed in higher dimensional flag spaces.

\begin{thm}  \label{theorem:rigidity}
\begin{itemize}
\item
Any local embedding of a CR structure into the flag space $\SL(3,\C)/B$ coincides locally with  $\phi_{CR}: {\mathbb S}^3 \rightarrow \SL(3,\C)/B$ up to a translation.  In particular, only spherical CR structures can be embedded.
\item Any local embedding of a path structure into the flag space $\SL(3,\C)/B$ coincides locally with  $\phi_{\R}: F=\SL(3,\R)/B_\R \rightarrow \SL(3,\C)/B$ up to a translation.  In particular, only flat path geometries can be embedded.
\end{itemize}
\end{thm}

\Pf  By the previous theorem the only CR structures which can be embedded are the spherical ones.  On the other hand, null curvature CR structures are known (\cite{C}) to be locally equivalent to ${\mathbb S}^3$ equipped to its standard structure.   By the theorem again, the null curvature structures admit embeddings which differ at most by a translation.  The proof in the case of path geometry is similar.
\EPf


\section{A global invariant}\label{section:global}


The second Chern class of the bundle $Y$ with connection form $\pi$ is given by
$$
c_2(Y,\pi)=\frac{1}{8\pi^2}\tr(\Pi\wedge \Pi).
$$

In the case of the connection form $\pi$ we obtain
$$
\left ( \begin{array}{ccc}
                        0 &   -\Phi^{2}      &  -\frac{1}{4}\Psi   \\

                         0  &  0    &  \frac{1}{2}\Phi^1  \\

                           0  &   0  &    0
               \end{array} \right )
               \wedge 
               \left ( \begin{array}{ccc}
                        0 &   -\Phi^{2}      &  -\frac{1}{4}\Psi   \\

                         0  &  0    &  \frac{1}{2}\Phi^1  \\

                           0  &   0  &    0
               \end{array} \right )
               =\left ( \begin{array}{ccc}
                        0 &   0     &  -\frac{1}{2}\Phi^1\wedge \Phi^2  \\

                         0  &  0    &  0\\

                           0  &   0  &    0
               \end{array} \right ).
$$
As $\Phi^1=Q^1\omega\wedge \omega^2$ and $\Phi^2=Q^2\omega\wedge \omega^1$ we have $\Pi\wedge \Pi=0$ and
$$
c_2(Y,\pi)=0.
$$

We include the proof of the next lemma although it is standard.

\begin{lem} The form 
$$
TC_2(\pi)=\frac{1}{8\pi^2}\left ( \tr(\pi\wedge \Pi)+\frac{1}{3}\tr(\pi\wedge \pi\wedge\pi)\right)=
\frac{1}{24\pi^2}\tr(\pi\wedge \pi\wedge\pi)
$$
is closed.

\end{lem}
\Pf 
Observe first that differentiating the curvature form we obtain $d\,\Pi=\Pi\wedge \pi-\pi\wedge \Pi$.  
Next we compute 
$$
d\,\tr (\Pi\wedge \pi)=\tr(d\,\Pi\wedge \pi+\Pi\wedge d\,\pi)
=\tr((\Pi\wedge \pi-\pi\wedge\Pi)\wedge \pi+\Pi\wedge (\Pi-\pi\wedge\pi))
$$
$$
=-\tr(\pi\wedge\Pi\wedge \pi)=0
$$
because
$$
\tr (\Pi\wedge \pi) =-\Phi^2\wedge \omega^1+\Phi^1\wedge \omega^2=0.
$$
Note that $\tr(\alpha\wedge \beta)=(-1)^{kl}\tr(\beta\wedge \alpha)$ if $\alpha$ and $\beta$ are two matrices of forms 
of degree $k$ and $l$ respectively.  Therefore, computing
$$
\frac{1}{3}d\,\tr(\pi\wedge \pi\wedge\pi)= \tr(d\,\pi\wedge \pi\wedge\pi)= \tr((\Pi-\pi\wedge\pi)\wedge \pi\wedge\pi)
$$
$$
=- \tr(\pi\wedge\pi\wedge \pi\wedge\pi)=0.
$$
\EPf

Remark that $0=c_2(Y,\pi)=d\, TC_2(\pi)$. 

\begin{dfn}\label{definition:global}
Suppose  that the fiber bundle $Y\rightarrow M$ is trivial and let $s: M\rightarrow Y$ be a section, we define then
$$
\mu=\int_M s^* TC_2(\pi)= \frac{1}{24\pi^2}\int_M s^*\tr(\pi\wedge \pi\wedge\pi).
$$
\end{dfn}

In principle that integral depends on the section but the following proposition shows that the integrand 
$$
s^* TC_2(\pi)
$$
defines
an element in the cohomology which does not depend on the section. 

 \begin{prop} 
 Suppose $s$ and $\tilde s$ are two sections.  Then
 $$
 \tilde s^*TC_2(\pi)-s^*TC_2(\pi)=-\frac{1}{8\pi^2}d\, s^*\tr (h^{-1}\pi\wedge d\,h ).
$$
 \end{prop}
 \Pf
 Fix the section $s$. Than there exists a map $h: M\rightarrow H$ such that
 $\tilde s = R_h\circ s$.  We have then
 $$
 \tilde s^*TC_2(\pi)=\frac{1}{24\pi^2}s^*\tr(R_h^*\pi\wedge R_h^*\pi\wedge R_h^*\pi).
 $$
From the formula
$$
{R_h}^*\pi= h^{-1}d\,h+ Ad_{h^{-1}}\pi,
$$
we obtain
$$
\tr(R_h^*\pi\wedge R_h^*\pi\wedge R_h^*\pi)=
$$
$$
\tr \left( h^{-1}d\,h\wedge  h^{-1}d\,h\wedge h^{-1}d\,h 
+ 3h^{-1}d\,h\wedge  h^{-1}\pi\wedge d\,h+3h^{-1}\pi\wedge \pi\wedge d\,h + \pi\wedge \pi\wedge\pi \right)
$$
$$
=\tr \left( - h^{-1}d\,h\wedge  d\,h^{-1}\wedge d\,h 
-3d\,h^{-1}\wedge \pi\wedge d\,h+3h^{-1}\pi\wedge \pi\wedge d\,h + \pi\wedge \pi\wedge\pi \right).
$$

\begin{lem} $\tr (h^{-1}d\,h\wedge  d\,h^{-1}\wedge d\,h )=0$.
\end{lem}

\Pf Observe that $ d\,h^{-1}\wedge d\,h$ is upper triangular with null diagonal.  
Moreover $h^{-1}d\,h$ is upper triangular and therefore the Lie algebra valued form also has zero diagonal.

\EPf

\begin{lem} $d\,\tr (h^{-1}\pi\wedge d\,h )=
\tr \left( d\,h^{-1}\wedge \pi\wedge d\,h-h^{-1}\pi\wedge \pi\wedge d\,h\right)$.
\end{lem}

\Pf
Compute
$d\tr (h^{-1}\pi\wedge d\,h )=\tr \left( d\,h^{-1}\wedge \pi\wedge d\,h+h^{-1}d\pi\wedge d\,h\right)$
$$
=\tr \left( d\,h^{-1}\wedge \pi\wedge d\,h+h^{-1}(\Pi -\pi\wedge \pi)\wedge d\,h\right)
$$

$$
=\tr \left( d\,h^{-1}\wedge \pi\wedge d\,h-h^{-1}\pi\wedge \pi\wedge d\,h\right)
$$
because $ \tr \left( h^{-1}\Pi \wedge d\,h\right)=0$ as in the previous lemma.

\EPf

The proposition follows from the two lemmas.

\EPf
\subsection{First variation}

We obtain in this section a first variation formula for the invariant $\mu$ when the flag structure 
is deformed through a smooth path.
Let $\mu(t)$ be the invariant defined as a function of the a parameter describing the deformation of the structure 
on a closed manifold $M$
and define $\delta\mu=\frac{d}{dt}\mu (0)$.

\begin{prop}$\label{proposition:variation}
\delta\mu= -\frac{1}{4\pi^2}\int_M s^*\tr(\dot \pi\wedge \Pi).
 $ 
 \end{prop}
 \Pf
 Differentiating $\mu(t)=\frac{1}{24\pi^2}\int_M s^*\tr(\pi_t\wedge \pi_t\wedge\pi_t)$ we have
 $$
 \delta\mu=\frac{1}{8\pi^2}\int_M s^*\tr(\pi\wedge \pi\wedge\dot{\pi}).
$$
Using the formula $\dot{\Pi}=d\dot{\pi}+\dot{\pi}\wedge\pi+\pi\wedge\dot{\pi}$ we write
$$
\tr({\pi}\wedge \dot{\Pi})=\tr(\pi\wedge d\dot{\pi} + 2\pi\wedge \pi\wedge\dot{\pi})
$$
and therefore
$$
\tr(\pi\wedge \pi\wedge\dot{\pi})=\frac{1}{2}\tr({\pi}\wedge \dot{\Pi}-\pi\wedge d\dot{\pi})
= \frac{1}{2}\tr(-\dot{\pi}\wedge {\Pi}-d\pi\wedge \dot{\pi}).
$$
Where, in the last equality, we used that on a closed manifold 
$ \int_M\tr({\pi}\wedge d\dot{\pi})=\int_M\tr(d{\pi}\wedge \dot{\pi})$ and that, differentiating
 $\tr({\pi}\wedge{\Pi})=0$, we have $\tr(\dot{\pi}\wedge {\Pi}+\pi\wedge d\dot{\Pi})=0$.
 
 Substituting $\Pi-\pi\wedge\pi= d\pi$ we obtain 
 $$
   \tr(\pi\wedge \pi\wedge\dot{\pi})= \frac{1}{2}\tr(-2\dot{\pi}\wedge {\Pi}-\pi\wedge\pi\wedge \dot{\pi})
$$
and therefore
$$
\delta\mu= \frac{1}{8\pi^2}\int_M s^*\tr(\pi\wedge \pi\wedge\dot{\pi})= 
-\frac{1}{4\pi^2}\int_M s^*\tr(\dot \pi\wedge \Pi).
$$

\EPf

Observe that an explicit computation gives 
$$
\tr(\dot \pi\wedge \Pi)= -\dot{\omega}^1\wedge \Phi^2 +\dot{\omega}^2\wedge \Phi^1 -\frac{1}{2}\dot{\omega}\wedge \Psi.
$$


\section{Pseudo flag geometry}\label{section:pseudo}

In this section we fix a contact form and obtain a reduction of the structure group of a path geometry.  We will obtain the relations 
between the invariants of the reduced structure to the original one.  This is similar to the reduction of a CR structure to a pseudo hermitian structure. 

We consider a  form $\theta$ on $T^\C$ such that $\ker \theta=T^1\oplus T^2$ is non-integrable.  
Define forms $Z^1$ and $Z^2$ on $T^\C$ satifying
$$
Z^1(T^1)\neq 0 \ Ê \  
{\mbox{and}}\ 
Z^2(T^2)\neq 0 ,
$$
$$
\ker \theta^1\supset T^2 \ Ê \  {\mbox{and}}\ 
 \  \ker \theta^2\supset T^1
$$
and such that $d\theta =Z^1\wedge Z^2$.

Fixing one choice, all others are given by $\theta^{1} = aZ^1$ and $\theta^{2} = \frac{1}{a}Z^2$, with $a\in \C^*$.
We consider now the $\C^*$ coframe bundle $X$ defined by the forms  $\theta^1,\theta^2,\theta$.
We have 
\begin{equation}\label{dtheta}
d\theta=\theta^1\wedge \theta^2.
\end{equation}

\begin{prop}\label{theta11}
There exist unique forms $\theta^1_1$, $\tau^1$ and $\tau^2$ on $X$ such that 
\begin{equation}\label{dtheta1}
d\theta^1=\theta^1\wedge\theta^1_1+\theta\wedge\tau^1
\end{equation}
\begin{equation}\label{dtheta2}
d\theta^2=-\theta^2\wedge \theta^1_1 +\theta\wedge \tau^2
\end{equation}
with $\theta^1_1=-\frac{da}{a}  \mod \theta^1,\ \theta^{2},\ \theta$ and $\tau^1\wedge \theta^2=\tau^2\wedge \theta^1=0$.
\end{prop}
\Pf
Define functions $z^i_{12},z^i_{j0}$ by
$$
dZ^i= z^i_{12}Z^1\wedge Z^2+ z^i_{10}Z^1\wedge \theta+ z^i_{20}Z^2\wedge \theta.
$$
Then 
$$
d\theta^1= \frac{da}{a}\wedge \theta^1+adZ^1= \frac{da}{a}\wedge \theta^1
+a\left (z^1_{12}Z^1\wedge Z^2+ z^1_{10}Z^1\wedge \theta+ z^1_{20}Z^2\wedge \theta\right).
$$
which can be written as
\begin{equation}\label{dtheta1}
d\theta^1=\theta^1\wedge \theta^1_1 +\theta\wedge \tau^1
\end{equation}
where $\theta^1_1=-\frac{da}{a} +z^1_{12}Z^2  -z^2_{12}Z^1$ (where we added a term in $Z^1$ in order to have a compatibility
with the formula for $d\theta^2$ bellow) and $\tau^1=-z^1_{10}\theta^1- z^1_{20}a^2\theta^2$.

Anagously, from

$$
d\theta^2= -\frac{da}{a}\wedge \theta^2+\frac{1}{a}dZ^2=- \frac{da}{a}\wedge \theta^2
+\frac{1}{a}\left (z^2_{12}Z^1\wedge Z^2+ z^2_{10}Z^1\wedge \theta+ z^2_{20}Z^2\wedge \theta\right)
$$
 we obtain
 \begin{equation}\label{dtheta2}
d\theta^2=-\theta^2\wedge \theta^1_1 +\theta\wedge \tau^2
\end{equation}
where  $\theta^1_1=-\frac{da}{a} +z^1_{12}Z^2-z^2_{12}Z^1 $ and $\tau^2=-z^2_{10}a^{-2}\theta^1- z^{2}_{20}\theta^2$.

Observe also that, differentiating equation \ref{dtheta} and using \ref{dtheta1} and \ref{dtheta2}, we obtain 
$\theta\wedge (\tau^1\wedge\theta^2-\tau^2\wedge \theta^1)=0$ which implies that, writing $\tau^i=\tau^i_1\theta^1+\tau^i_2\theta^2$,
$$
\tau^1_1+\tau^2_2=0.
$$
Now, if $\theta'^1_1$, $\tau'^1$ and $\tau'^2$ are other forms satisfying the equations, then from the above equations
we obtain
$$
\theta^1_1-\theta'^1_1=A\theta \ \ \mbox{and}\ \  \tau^1_1-\tau'^1_1=A.
$$
 Choosing an appropriate $A$ we can therefore fix $\tau^1_1=0$ and the forms $\theta^1_1$, $\tau^1$ and $\tau^2$ are uniquely determined as claimed.
 \EPf
 
 Differentiating equations \ref{dtheta1} and \ref{dtheta2} we obtain
$$
\theta^1\wedge d\theta^1_1+\theta\wedge (d\tau^1-\tau^1\wedge \theta^1_1)=0
$$
$$
-\theta^2\wedge d\theta^1_1+\theta\wedge (d\tau^2+\tau^2\wedge \theta^1_1)=0
$$
and therefore
\begin{equation}\label{dtheta11}
d\theta^1_1=R\theta^1\wedge \theta^2+W^1\theta^1\wedge\theta+W^2\theta^2\wedge\theta
\end{equation}
\begin{equation}\label{dtau1}
d\tau^1-\tau^1\wedge \theta^1_1=-W^2\theta^1\wedge \theta^2+S^1_1\theta\wedge\theta^1+S^1_2\theta\wedge\theta^2
\end{equation}
\begin{equation}\label{dtau2}
d\tau^2+\tau^2\wedge \theta^1_1=-W^1\theta^1\wedge \theta^2+S^2_1\theta\wedge\theta^1+S^2_2\theta\wedge\theta^2
\end{equation}
Moreover, differentiating equation $\tau^1\wedge \theta^2=0$ and  $\tau^2\wedge \theta^1=0$ 
we obtain
$$
S^1_1=S^2_2=\tau^1_2\tau^2_1.
$$

\subsection{Curvature identities}

Differentiating equation \ref{dtheta11} one gets
$$
dR\wedge  \theta^1\wedge \theta^2+ (dW^1-W^1\theta^1_1)\wedge \theta^1\wedge \theta
 (dW^2+W^2\theta^1_1)\wedge \theta^2\wedge \theta=0.
 $$ 
 Writing
 $$
 dR=R_0\theta +R_1\theta^1+R_2\theta^2,
 $$
$$
 dW^1-W^1\theta^1_1=W^1_0\theta+W^1_1\theta^1+W^1_2\theta^2
 $$
 and
 $$
  dW^2+W^2\theta^1_1=W^2_0\theta+W^2_1\theta^1+W^2_2\theta^2
 $$
 Then
 $$
 R_0=W^1_2-W^2_1.
 $$
 Differentiating equation \ref{dtau2} and writing $dR_0=R_{00}\theta+R_{01}\theta^1+R_{02}\theta^2$
, one gets
$$
dR_1-R_1\theta^1_1+R_2\tau^2_1\theta-\frac{1}{2}R_0\theta^2=R_{01}\theta+R_{11}\theta^1+R_{12}\theta^2
 $$
 and
 $$
dR_2+R_2\theta^1_1+R_1\tau^1_2\theta+\frac{1}{2}R_0\theta^1=R_{02}\theta+R_{12}\theta^1+R_{22}\theta^2
 $$
 We also obtain differentiating \ref{dtau1} and \ref{dtau2}
 $$
 d\tau^1_2+2\tau^1_2\theta^1_1=-W^2\theta^1 +S^1_2\theta \mod \theta^2
 $$
 and
 $$
 d\tau^2_1-2\tau^2_1\theta^1_1=W^1\theta^2 +S^2_1\theta \mod \theta^1.
 $$
 \subsection{Embedding $X\rightarrow Y$}\label{section:embeddding}
 
 Recall that $X$ is the coframe bundle of forms $(\theta,\theta^1,\theta^2)$  over $M$. We chose a section of this bundle.  The forms over $M$ will also be denoted by $(\theta,\theta^1,\theta^2)$.
 The goal now is to obtain an immersion $X\rightarrow Y$.  
 Let $s:M\rightarrow Y$ such that $s^* \omega=\theta$ and write
 $$
 s^*\phi= A_1\theta^1+A_2\theta^2+A_0\theta.
 $$
 for functions $A_i$ on $M$. To choose a section we will impose
 that $s^*\phi=0$.  For that sake we start with a particular section and move it using the action of the structure group $H$.
 
Consider$$
h=
\left( \begin{array}{ccc}

                        a  	&    	c      	&       e \\

                       0				 &     \frac{1}{ab}		&		f\\
                       0   				&  	0 				&		b
                      
       \end{array} \right ) 
$$
which gives
\begin{equation}\label{hdh}
h^{-1}d\,h=
					\left( \begin{array}{ccc}          
          a^{-1}d\,a  	&    	a^{-1}d\, c+c(\frac{d\,b}{ab}+\frac{d\,a}{a^2}) 	&      a^{-1}(d\,e)-bc(d\,f)+(cf-\frac{e}{ab})d\, b	\\
           0				 	&     -\frac{d\,(ab)}{ab}		&		ab(d\,f)-af(d\,b)															\\
           0   				&  	0 						&		\frac{d\, b}{b}
                       \end{array} \right ).
\end{equation}
We will use the formula
$$
{R_h}^*\pi= h^{-1}d\,h+ Ad_{h^{-1}}\pi.
$$

If $h: M\rightarrow H$ is given by
$$
\left( \begin{array}{ccc}

                        1 	&    	c      	&       e \\

                       0				 &     1		&		f\\
                       0   				&  	0 				&		1
                      
       \end{array} \right ) 
$$ 
then, using formula \ref{hdh} and formula \ref{adjunta} for $Ad_{h^{-1}}\pi$ in the expression of ${R_h}^*\pi$ we obtain 

$$
\tilde \phi =\phi +c\omega^1+2f\omega^2 +(4e-2cf)\omega.
$$
Observe that $a=b$ is imposed by the condition  $s^* \omega=\theta$.
Starting with any section defining functions $A_1$, $A_2$ and $A_3$ we obtain a new section 
by acting by a section of $h: M\rightarrow H$ given by
 $c=-A_1,\ f=-A_2/2,\ e=(A_1A_2-A_0)/4$.
 In that case we have
  $$
 s^*\tilde\phi=0.
 $$
 
 The maps from $E$ to the fiber group of $Y\rightarrow E$, such that $s^*{R_h}^*\tilde \phi=0$ are fixed, because
 ${R_h}^*\tilde\phi=\tilde\phi+c\,\omega^1+2{f}{}\,\omega^2+({4e}-2cf)\,\omega$,
so $c=f=e=0$.

After fixing $c,f,e$,  we allow  maps $h:M\rightarrow H$ acting by $R_h$ on $Y\rightarrow M$ by 
elements of the form
%
 $$
h=\left( \begin{array}{ccc}

                        a  	&    	0      	&       0 \\

                       0				 &     \frac{1}{a^2}		&		0\\
                       0   				&  	0 				&		a
                      
       \end{array} \right ) 
$$
so that the form $\theta$ be preserved.  This gives the embedding of $X$ into $Y$.

We may suppose that $s^*\phi=0$ and then obtain the following equations by pulling back to $M$ the structure equations
on $Y$:
\begin{equation}
d\theta=\theta^1\wedge \theta^2
\end{equation}
\begin{equation}\label{dtheta1'}
d\theta^1=\theta^1\wedge \omega^1_1+\theta\wedge \phi^1
\end{equation}
 \begin{equation}\label{dtheta2'}
d\theta^2=-\theta^2\wedge \omega^1_1+\theta\wedge \phi^2
\end{equation}
\begin{equation}\label{psi}
\theta^1\wedge \phi^2-\theta^2\wedge \phi^1+\theta\wedge \psi=0
\end{equation}
\begin{equation}\label{domega11'}
d\omega^1_1-\frac{3}{2}\theta^1\wedge \phi^2-\frac{3}{2}\theta^2\wedge \phi^1=0
\end{equation}
\begin{equation}\label{dphi1'}
d\phi^1-\phi^1\wedge \omega^1_1-\frac{1}{2}\theta^1\wedge \psi=Q^1\theta\wedge \theta^2
\end{equation}
\begin{equation}\label{dphi2'}
d\phi^2+\phi^2\wedge \omega^1_1-\frac{1}{2}\theta^2\wedge \psi=Q^2\theta\wedge \theta^1
\end{equation}
\begin{equation}
d\psi-2\phi^1\wedge \phi^2=U_1\theta^1\wedge \theta+U_2\theta^2\wedge \theta
\end{equation}

In the formulae above we write the pull back of any form $\alpha$ defined on $Y$ using the same notation $\alpha$. 
It follows from equations \ref{dtheta1'} and \ref{dtheta2'}, comparing with proposition \ref{theta11} that
$$
\omega^1_1=\theta^1_1+c\theta
$$
and therefore
\begin{equation}
d\theta^1=\theta^1\wedge \theta^1_1+\theta\wedge( \phi^1-c\theta^1)
\end{equation}
 \begin{equation}
d\theta^2=-\theta^2\wedge \theta^1_1+\theta\wedge (\phi^2+c\theta^2)
\end{equation}
which implies
$$
\phi^1=c\theta^1+E^1\theta +\tau^1,
$$
$$
\phi^2=-c\theta^2+E^2\theta +\tau^2,
$$
Substituting the formulae above in equation \ref{psi}, and using again proposition \ref{theta11} we obtain
$$
\theta\wedge (\psi-E^2\theta^1+E^1\theta^2)=0
$$
and therefore
$$
\psi= E^2\theta^1-E^1\theta^2+G\theta.
$$
Substituting the expressions of $\omega^1_1, \phi^1$ and $\phi^2$ in equation \ref{domega11'} 
and using equation \ref{dtheta11} we obtain
$$
(R+4c)\,\theta^1\wedge\theta^2 + (W^1-\frac{3}{2}E^2)\,\theta^1\wedge \theta +(W^2-\frac{3}{2}E^1)\,\theta^2\wedge \theta+dc\wedge \theta=0.
$$
This implies
$$
c=-\frac{R}{4}
$$
writing, as before, $dR=R_0\theta +R_1\theta^1+R_2\theta^2$ and substituting in the above expression we get
$$
 (W^1-\frac{3}{2}E^2-\frac{1}{4}R_1)\,\theta^1\wedge \theta +(W^2-\frac{3}{2}E^1-\frac{1}{4}R_2)\,\theta^2\wedge \theta=0
$$
which implies
$$
E^2=\frac{2}{3}(W^1-\frac{1}{4}R_1)
$$
and
$$
E^1=\frac{2}{3}(W^2-\frac{1}{4}R_2).
$$
We now write equation \ref{dphi1'}, substituting the expressions for $\phi^1$ and $\omega^1_1$ obtained above:
$$
0= d\tau^1-\frac{1}{4}(R_0\theta+R_2\theta^2)\wedge\theta^1
-\frac{1}{4}R(\theta^1\wedge \theta^1_1+\theta\wedge \tau^1)+
d E^1\wedge \theta +E^1\theta^1\wedge \theta^2
$$
$$
-(\tau^1-\frac{1}{4}\theta^1+E^1\theta)\wedge (\theta^1_1-\frac{1}{4}R\theta)
-\frac{1}{2}\theta^1\wedge (E^2\theta^1-E^1\theta^2+G\theta)-Q^1\theta\wedge\theta^2
$$
Using curvature identities in order to simplify the above expression we obtain after
a computation
$$
 0=\theta\wedge \theta^1(S^1_1-\frac{1}{3}R_0+\frac{1}{16}R^2+\frac{1}{2}G-\frac{2}{3}W_1^2+\frac{1}{6}R_{21})
 +\theta\wedge \theta^2(S1_2-\frac{1}{2}R\tau^1_2-Q^1-\frac{2}{3}W^2_2+\frac{1}{6}R_{22})
 $$
 Therefore
 $$
 S^1_1-\frac{1}{3}R_0+\frac{1}{16}R^2+\frac{1}{2}G-\frac{2}{3}W_1^2+\frac{1}{6}R_{21}=0
 $$
 and
 $$
 S^1_2-\frac{1}{2}R\tau^1_2-Q^1-\frac{2}{3}W^2_2+\frac{1}{6}R_{22}=0.
 $$
Analogously, from equation \ref{dphi2'}, we obtain the identities
$$
 S^2_2+\frac{1}{3}R_0+\frac{1}{16}R^2+\frac{1}{2}G-\frac{2}{3}W_2^1+\frac{1}{6}R_{12}=0
 $$
 and
 $$
 S^2_1+\frac{1}{2}R\tau^2_1-Q^2-\frac{2}{3}W^1_1+\frac{1}{6}R_{11}=0.
 $$

 \subsubsection{The global invariant}
 
 A simple computation gives the formula
 $$
 \tr (\pi\wedge\pi\wedge \pi)=\frac{3}{2}(\omega\wedge \phi +\omega^1\wedge \omega^2)\wedge \psi
+3\omega\wedge \phi^1\wedge \phi^2+3\omega^1\wedge (\frac{1}{2}\phi+\omega^1_1)\wedge\phi^2
-3\omega^2\wedge (\frac{1}{2}\phi-\omega^1_1)\wedge\phi^1.
$$
Therefore using the embedding of the previous section we obtain by a computation:
\begin{prop} \label{proposition:invariant}
For a section $s : M\rightarrow Y$ factoring through an embedding   of $X$ into $Y$
such as $s^*\phi=0$ we have
$$
s^*TC_2(\pi)=\frac{1}{8\pi^2}\left( \theta\wedge \theta^1\wedge \theta^2(\frac{1}{2} G+\frac{1}{16}R^2
-\tau^1_2\tau^2_1)
+\theta^1_1\wedge(E^2\theta\wedge\theta^1+
E^1\theta\wedge\theta^2-\frac{R}{2} \theta^1\wedge\theta^2)\right).
$$
\end{prop}

\section{Homogeneous flag structures on $\SU(2)$ }

Let $\alpha, \beta, \gamma$ be a basis of left invariant 1-forms defined
on $\SU(2)$ with
$$
d\alpha= -\beta\wedge \gamma,\ \ \ d\beta= -\gamma\wedge \alpha,\ \ \ d\gamma= -\alpha\wedge \beta\ \ \ 
$$
We define a pseudo flag structure choosing a map from $\SU(2)$ to $\SL(2,\C)$:
$$
\theta=\gamma, \ \ \ Z^1=r_1\beta+r_2\alpha, \ \ \ Z^2=s_1\beta+s_2\alpha,
$$
with $r_1s_2-r_2s_1=1$.
Then
$$
d\theta=Z^1\wedge Z^2.
$$
 
In the case the map $\SU(2)\rightarrow \SL(2,\C)$ is constant, from 
$\beta =s_2Z^1-r_2Z^2$ and $\alpha =-s_1Z^1+r_1Z^2$,
we obtain
$$
dZ^1=r_1d\beta +r_2 d\alpha=\theta \wedge \left( xZ^1+yZ^2\right)
$$
and analogously,
$$
dZ^2=\theta \wedge \left( zZ^1-xZ^2\right),
$$
 where
$$
x=r_1s_1+r_2s_2,\ \ \ y=-(r_1^2+r_2^2), \ \ \ z=s_1^2+s_2^2.
$$
Observe that $x^2+yz=-1$.
Then for a pseudo flag structure with coframes obtained from the tautological forms $\theta^1=aZ^1, \theta^2=a^{-1}Z^2$
$$
d\theta^1=\theta^1\wedge (-\frac{da}{a}-x\theta) + \theta\wedge (y a^2 \theta^2)
$$
$$
d\theta^2=-\theta^2\wedge (-\frac{da}{a}-x\theta) + \theta\wedge (z a^{-2} \theta^1)
$$
From Proposition \ref{theta11} we have
$$
\theta^1_1=-\frac{da}{a}-x\theta,
$$
$$
\tau^1=y a^2 \theta^2, \ \ \ \tau^2=z a^{-2} \theta^1.
$$
and therefore
$$
d\theta^1_1=-x\theta^1\wedge \theta^2
$$
so that $R=-x, W^1=W^2=0$.
In order to compute the curvature invariants from the pseudo flag structure we 
use the embedding in section \ref{section:embeddding}.
We compute first
$$
d\tau^1=ya^2\theta^2\wedge\theta^1_1+yz\theta\wedge \theta^1-2a^2xy\theta\wedge \theta^2
$$
and
$$
d\tau^2=-za^{-2}\theta^1\wedge\theta^1_1+2xza^{-2}\theta\wedge \theta^1+yz\theta\wedge \theta^2.
$$

Now, as $d\tau^1-\tau^1\wedge \theta^1_1=-W^2\theta^1\wedge \theta^2+S^1_1\theta\wedge\theta^1+S^1_2\theta\wedge\theta^2$
and $d\tau^2+\tau^2\wedge \theta^1_1=-W^1\theta^1\wedge \theta^2+S^2_1\theta\wedge\theta^1+S^2_2\theta\wedge\theta^2$
(cf. \ref{dtau1},\ref{dtau2}) we obtain
$$
S^1_1=yz, \ \ S^1_2=-2a^2xy, \ \ S^2_1=2xza^{-2}, \ \ S^2_2=yz.
$$
From the embedding equations we have
$$
Q^1=S^1_2-\frac{1}{2}R\tau^1_2-\frac{2}{3}W^2_2+\frac{1}{6}R_{22}= -\frac{3}{2}xya^2
$$
and, analogously
$$
Q^2=S^2_1+\frac{1}{2}R\tau^2_1-\frac{2}{3}W^1_1+\frac{1}{6}R_{11}= \frac{3}{2}xza^{-2}.
$$
We conclude that $Q^1=Q^2=0$ if and only if $x=0$ or $y=z=0$.

\subsubsection{The global invariant}
We compute, using \ref{proposition:invariant}, the global invariant for the family of structures defined on $\SU(2)$
We have $R=-x$, $\tau^1_2=a^2y$,  $\tau^2_1=a^{-2}z$, $E^1=E^2=0$ and $G=-2yz-\frac{1}{8}x^2$.
Therefore for a section $s:\SU(2)\rightarrow Y$ as above we obtain
$$
s^*TC_2(\pi)=-\frac{1}{8\pi^2}\gamma\wedge \beta\wedge \alpha(2yz+\frac{1}{2}x^2).
$$



\begin{flushleft}
  \textsc{E. Falbel\\
  Institut de Math\'ematiques \\
  de Jussieu-Paris Rive Gauche \\
CNRS UMR 7586 and INRIA EPI-OURAGAN \\
 Sorbonne Universit\'e, Facult\'e des Sciences \\
4, place Jussieu 75252 Paris Cedex 05, France \\}
 \verb|elisha.falbel@imj-prg.fr|
 \end{flushleft}
\begin{flushleft}
  \textsc{J. M.  Veloso\\
  Faculdade de Matem\' atica - ICEN\\
Universidade Federal do Par\'a\\66059 - Bel\' em- PA - Brazil}\\
  \verb|veloso@ufpa.br|
\end{flushleft}


\begin{thebibliography}{ZZ99}

\bibitem[Ba]{Ba} Barbot, T. :
Flag structures on Seifert manifolds.
Geom. Topol. 5 (2001), 227-266. 

\bibitem[BFG]{BFG} Bergeron, N. , Falbel, E. , Guilloux, A. : Tetrahedra of flags, volume and homology of $\SL(3)$. Geom. Topol. 18 (2014), no. 4, 1911-1971.

\bibitem[BHR]{BHR} Biquard, O. , Herzlich, M. , Rumin, M. : 
Diabatic limit, eta invariants and Cauchy-Riemann manifolds of dimension 3. 
Ann. Sci. \'Ecole Norm. Sup. (4) 40 (2007), no. 4, 589-631.

\bibitem[Bo]{Bo} Borrelli, V. : On totally real isotopy classes. 
Int. Math. Res. Not. 2002, no. 2, 89-109. 

\bibitem[B]{B} Bryant, R. :  
\'Elie Cartan and geometric duality. Preprint 1998.

\bibitem[BGH]{BGH} Bryant, R. , Griffiths, P. , Hsu, L. :
Toward a geometry of differential equations. Geometry, topology, and physics, 1-76, 
Conf. Proc. Lecture Notes Geom. Topology, IV, Int. Press, Cambridge, MA, 1995. 

 
\bibitem[BE]{BE} Burns, D. , Epstein, C. L. A global invariant for three-dimensional CR-manifolds. Invent. Math. 92 (1988), no. 2, 333-348.

\bibitem[BS]{BS} Burns, D. , Shnider, S. : Real Hypersurfaces in complex manifolds.  Proceedings of
Symposia in Pure Mathematics, 30  (1977), 141-168.

\bibitem[BS1]{BS1} D. Burns , S. Shnider : Spherical hypersurfaces in complex manifolds. Invent. Math. 33 (1976), 223Ð246


\bibitem[Ca]{Ca} Cartan, E. : Sur les vari\'et\'es ˆ connexion projective. Bull. Soc. Math. France 52 (1924), 205-241.

\bibitem[C]{C} Cartan, E. : Sur la g\' eom\' etrie pseudo-conforme des
 hypersurfaces de deux variables complexes, I.  Ann. Math. Pura Appl., (4) 11 (1932) 17-90
 (or Ouevres II, 2, 1231-1304); II, Ann. Scuola Norm. Sup. Pisa, (2) 1 (1932) 333-354
 (or Ouevres III, 2, 1217-1238).
 
 \bibitem[CL]{CL} Cheng, J. H. ,  Lee, J. M. : The Burns-Epstein invariant and deformation of CR structures. Duke Math. J. 60 (1990), no. 1, 221-254.

\bibitem[CM]{CM} Chern, S. S. , Moser, J. :
 Real Hypersurfaces in Complex Manifolds. Acta Math. 133 (1974) 219-271.
 
\bibitem[DF]{DF} Deraux, M. , Falbel, E. : Complex hyperbolic geometry of the figure-eight knot.  
Geom. Topol. 19 (2015), no. 1, 237-293.
 
\bibitem[FS]{FS} Falbel, E. , Santos Thebaldi, R. : A flag structure on a cusped hyperbolic 3-manifold. Pacific J. Math. 278 (2015), no. 1, 51-78.

  \bibitem[Fo]{Fo}
  Forstneric, F. : On totally real embeddings into $\C^n$. Exposition. Math. 4 (1986), no. 3, 243-255.
  
 \bibitem[IL]{IL}  Ivey, T. A. , Landsberg, J. M. : Cartan for beginners: differential geometry via moving frames and exterior differential systems. Graduate Studies in Mathematics, 61. American Mathematical Society, Providence, RI, 2003.


\bibitem[J]{J} Jacobowitz, H. : An Introduction to CR Structures. vol. 32, Mathematical Surveys and Monographs, American Mathematical Society, Providence, Rhode Island, 1990.

\bibitem[S]{S} Schwartz, R. E. : Spherical CR geometry and Dehn surgery.  Ann. of Math. Stud., 165, Princeton University Press, Princeton, NJ, 2007.

\bibitem[W]{W} Webster, S. M. : Pseudo-Hermitian structures on a real hypersurface. J. Differential Geom. 13 (1978), no. 1, 25-41.

\end{thebibliography}
\end{document}